\documentclass{article}
\usepackage{amsmath,amsfonts,amsthm}

\theoremstyle{plain}
\newtheorem{theorem}{Theorem}[section]
\newtheorem{corollary}{Corollary}
\newtheorem{lemma}{Lemma}

\theoremstyle{definition}
\newtheorem{example}{Example}
\newtheorem{definition}{Definition}

{%
\begin{enumerate}}%
{\end{enumerate}%
}%

\newcommand{\C}{\text{$\mathbb C$}}
\newcommand{\N}{\text{$\mathbb N$}}

\renewcommand{\frak}[1]{\text{$\mathfrak{#1}$}}
\newcommand{\J}{\text{$\mathcal{J}$}}

\newcommand{\ga}{\text{$\alpha$}}

\newcommand{\gO}{\text{$\Omega$}}
\newcommand{\gf}{\text{$\varphi$}}

\newcommand{\del}{\text{$\partial$}}

\newcommand{\delbar}{\text{$\overline{\partial}$}}
\newcommand{\tensor}{\otimes}
\newcommand{\im}{\mathrm{Im}\,}

\newcommand{\Ker}{\mathrm{Ker}\,}
\newcommand{\mc}[1]{\text{$\mathcal{#1}$}}

\newcommand{\into}{\longrightarrow}
\newcommand{\noqed}{\let\qed\relax}

\newcommand{\Gg}{\mathfrak{g}}

\newcommand{\IP}[1]{\langle #1 \rangle}

\newcommand{\gcs}{generalized complex structure}

\newcommand{\gcss}{generalized complex structures}

\newcommand{\gcm}{generalized complex manifold}
\newcommand{\gcms}{generalized complex manifolds}

\newcommand{\gk}{generalized K\"ahler}
\newcommand{\gks}{generalized K\"ahler structure}

\newcommand{\wrt}{with respect to}

\newcommand{\Cour}[1]{[\![#1]\!]}
\newcommand{\Jac}{\mathrm{Jac}\,}
\newcommand{\Cliff}{\mathrm{Cliff}\,}

\begin{document}

\title{Formality in generalized K\"ahler geometry}
\author{Gil R. Cavalcanti\thanks{{\tt gil.cavalcanti@maths.ox.ac.uk}} \\
        Mathematical Institute\\
     St. Giles 24-29\\
        Oxford, OX1 3LB, UK\\
[0.1cm]}

\maketitle

\abstract{We prove that no nilpotent Lie algebra admits an invariant generalized K\"ahler structure. This is done by showing that a certain differential graded algebra associated to a generalized complex manifold is formal in the generalized K\"ahler case, while it is never formal for a generalized complex structure on a nilpotent Lie algebra.}

\section*{Introduction}

Generalized K\"ahler manifolds, as introduced by Gualtieri \cite{Gu03}, have received recently a lot of attention from both physicists and mathematicians. From the physics point of view, they are the general  solution to the $(2,2)$ supersymmetric sigma model \cite{GHR84}, while for mathematicians they represent interesting examples of bihermitian structures \cite{Gu03}. A classification of manifolds which admit such structures was achieved in four dimensions \cite{ApGu06} and finding concrete examples of such structures has been a driving force in this area \cite{BCG05,Hi06,Hi05,LiTo05}.

However very little is known about the differential topology of \gk\ manifolds. This is  despite of the fact that their better known relatives, K\"ahler manifolds, have very restrictive topological properties, e. g., they have even `odd Betti numbers', satisfy the strong Lefschetz property and are   formal in the sense of Sullivan \cite{Su78,DGMS75}. The last property, for example, can be easily used to prove that no nilpotent Lie algebra has a K\"ahler structure \cite{Has89}.  The point of this paper is to  prove the analogous result for \gks s. We achieve this goal by providing a formality result for  \gk\ manifolds.

Generalized K\"ahler manifolds are a special type of \gcms, and as  such many of their properties stem from general properties of \gcss. For example, every \gcs\ induces a decomposition of the  forms analogous to the $(p,q)$-decomposition of forms on a complex manifold,  and this decomposition causes the exterior derivative to decompose as $d = \del +\delbar$. By studying Hodge theory on \gk\ manifolds, Gualtieri showed  in \cite{Gu04} that in a \gk\ manifold a series of  $\del\delbar$-lemma-like hold. 

Given that formality of K\"ahler manifolds is a consequence of the $\del\delbar$-lemma, one might expect that Gualtieri's result implies formality of \gk\ manifolds by the same argument from \cite{DGMS75}. However, in a \gcm, the operators $\del$ and $\delbar$ are not derivations (an important fact in for the formality theorem for K\"ahler manifolds) and not only does that spoil the proof, but there are examples of nonformal \gcms\ which {\it satisfy} the $\del\delbar$-lemma \cite{Ca04}.

In this paper we advance on this problem by abandoning the differential algebra $(\Omega^{\bullet}(M),d)$, and hence the question of whether $M$ is formal, and looking at a different DGA.

A \gcs\ can be described in terms of a Lie algebroid $L \subset (TM \oplus T^*M)\tensor \C$ and hence the exterior algebra $(\Gamma(\wedge^{\bullet}L^*),d_L)$ is a DGA. The key observation for our result is that on a \gcm\ with holomorphically trivial canonical bundle, this algebra is isomorphic, as a differential complex, to $(\Omega_{\C}^{\bullet}(M),\delbar)$.  In a \gk\ manifold, the operator $\delbar$ decomposes further $\delbar = \delta_+ +\delta_-$, giving rise to a decomposition $d_L = \del_L + \delbar_L$, with the advantage that, unlike $\del$, $\delbar$, $\delta_+$ or $\delta_-$, the operators  $\del_L$  and $\delbar_L$ are derivations. Hence, in this setting, using the correspondence between the different operators and Gualtieri's $\del\delbar$-lemmas we can prove formality of  $(\Gamma(\wedge^{\bullet}L^*),d_L)$.

This result provides concrete differential-topological obstructions and  allows us to prove that there are no invariant \gks s on nilpotent Lie algebras.

I would like to thank Marco Gualtieri  for his suggestions and inspired advice as well as Marisa Fern\'andez, Nigel Hitchin and Simon Salamon for suggestions on the first manuscript. I also thank IMPA for their hospitality while writing this paper. This research is supported by EPSRC.

\section{Differential graded algebras}

In this section we give a lightening review formality  for differential graded algebras and recover the well known fact that the DGA associated to a nontrivial nilpotent Lie algebra is not formal.

\begin{definition}
A {\it differential graded algebra}, or {\it DGA} for short, is an $\N$ graded vector space $\mc{A}^{\bullet}$, endowed with a product and a differential $d$ satisfying:
\begin{enumerate}
\item The product maps $\mc{A}^{i}\times \mc{A}^{j}$ to $\mc{A}^{i+j}$ and is graded commutative:
$$a \cdot b = (-1)^{ij} b\cdot a;$$
\item The differential has degree 1, $d:\mc{A}^{k} \into \mc{A}^{k+1}$, and squares to zero;
\item The differential is a derivation: for $a \in \mc{A}^{i}$ and $b \in \mc{A}^{j}$
$$d(a\cdot b) = da \cdot b + (-1)^i a \cdot db.$$
\end{enumerate}
\end{definition}

The cohomology of a DGA is defined in the standard way and naturally inherits a grading and a product, making it into a DGA on its own with $d=0$. A {\it morphism} of differential graded algebras is a map preserving the structure above, i.e., degree, products and differentials. Any morphism of DGAs $\gf:\mc{A} \into \mc{B}$ gives rise to a morphism of cohomology $\gf^*:H^{\bullet}(\mc{A}) \into H^{\bullet}(\mc{B})$.
A morphism $\gf$ is a {\it quasi-isomorphism} if the induced map in cohomology is an isomorphism.

Given a DGA, $\mc{A}$, for which $H^k(\mc{A})$ is finite dimensional for every $k$, one can construct another differential  graded algebra that captures all the information about the differential and which is minimal in the following sense.

\begin{definition}
A DGA $(\mc{M},d)$ is {\it minimal} if it is free as a DGA (i.e. polynomial in even degree and skew symmetric in odd degree) and has generators $e_1, e_2 \dots, e_n, \dots$ such that
\begin{enumerate}
\item The degree of the generators form a weakly increasing sequence of positive numbers;
\item There are finitely many generators in each degree;
\item The differential satisfies $de_i \in \wedge\mbox{span}\{e_1, \dots, e_{i-1}\}$.
\end{enumerate}
A {\it minimal model} for a differential graded algebra $\mc{A}$ is a minimal DGA, $\mc{M}$, together with a quasi-isomorphism $\psi: \mc{M} \into \mc{A}$.
\end{definition}

Since the cohomology of a DGA is also a DGA we can also construct its minimal model. The minimal models  for $\mc{A}$ and $H^{\bullet}(\mc{A})$ are not the same in general.

\begin{definition}
A DGA $\mc{A}$ is {\it formal} if it has the same minimal model as its cohomology, or equivalently, there is a quasi-isomorphism $\psi:\mc{M} \into H^{\bullet}(\mc{A})$, where $\mc{M}$ is the minimal model of $\mc{A}$. A manifold $M$ is {\it formal} if $(\Omega^{\bullet}(M),d)$ is formal.
\end{definition}

\begin{example}[Nilpotent Lie algebras \cite{Has89}]\label{ex:nilpotent lie algebras} A typical example of nonformal DGA can be obtained from a finite dimensional nilpotent Lie algebra $\frak{g}$ with nontrivial bracket. The Lie bracket induces a differential $d$ on  $\wedge^{\bullet} \Gg^*$ making it into a DGA. Furthermore, $\Gg^*$ is filtered by $\Gg^*_1 = \ker d$ and 
$$\Gg^*_i = \{v \in \Gg^*: dv \in \wedge^2 \Gg^*_{i-1}\}.$$
Nilpotency implies that $\Gg^*_s = \Gg^*$ for some $s$.
Let $\{e^1, \cdots, e^n\}$ be a basis for $\Gg^*$ compatible with this filtration. Then
$$de^i \in \wedge^2\mathrm{span}\{e^1, \cdots, e^{i-1}\}.$$
showing that $(\wedge^{\bullet}\Gg^*,d)$ is minimal.

Since the bracket is nontrivial, $de^n \neq 0$ and hence one can see that $e^1 \wedge \cdots \wedge e^{n-1}$ is exact and $e^1 \wedge \cdots \wedge e^{n}$ is a volume element and therefore represents a nontrivial cohomology class. If $(\wedge\Gg^*,d)$ was formal, there would be a map $\psi: (\wedge^{\bullet}\Gg^*,d) \into H^{\bullet}(\Gg)$, but
$$0 \neq  \psi(e^1 \wedge \cdots \wedge e^n) =  \psi(e^1 \wedge \cdots \wedge e^{n-1}) \cdot \psi(e^n) = 0 \cdot \psi(e^n) = 0.$$
So there is no such $\psi$ and $\wedge^{\bullet} \Gg^*$ is not formal.
\end{example}

\section{Generalized complex structures and Lie algebroids}

In this section we recall the definition of 
generalized complex structures and their relation to Lie algebroids, following~\cite{Gu03}.

Given a closed 3-form $H$ on a manifold $M$, we define the {\it
Courant bracket} of sections of  $T\oplus T^*$,  the sum of the tangent
and cotangent bundles, by
$$\Cour{X+\xi,Y+\eta}= [X,Y] + \mc{L}_X \eta -\mc{L}_Y\xi -\tfrac{1}{2}d(\eta(X) - \xi(Y)) + i_Y i_X H.$$
The bundle $T\oplus T^*$ is also endowed with a natural symmetric
pairing of signature $(n,n)$:
$$ \IP{X+\xi,Y+\eta} = \frac{1}{2} (\eta(X) + \xi(Y)).$$

\begin{definition}
A {\it \gcs} on a manifold with closed 3-form  $(M,H)$ is a complex
structure on the bundle $T\oplus T^*$ which preserves the natural
pairing and whose $i$-eigenspace is closed under the Courant
bracket.
\end{definition}

A \gcs\ can be fully described in terms of its $i$-eigenspace $L$,
which is a maximal isotropic subspace of $T_{\C} \oplus T^*_{\C}$
satisfying $L \cap \overline{L} = \{0\}$.

Two extreme examples of \gcss, with $H=0$, are given by complex and symplectic structures: in a complex manifold we let  $L = T^{0,1}\oplus T^{*1,0}$ and in a symplectic manifold we let $L = \{X -i \omega(X):X \in T_{\C}\}$, where $\omega$ is the symplectic form. What distinguishes these structures is their {\it type} which is the dimension of the kernel of  $\pi:L {\into} T_{\C}$. So, a complex structure on $M^n$  has type $n$ at all points and symplectic structures have type zero at all points.

The Courant bracket does not satisfy the Jacobi identity. Instead we have the relation for the Jacobiator
$$\Jac(A,B,C):= \Cour{\Cour{A,B},C} + c.p.  = \tfrac{1}{3}d(\IP{\Cour{A,B},C} +c.p.), $$
where $c.p.$ stands for cyclic permutations. However, the identity above also shows that the Courant bracket induces a Lie bracket  when restricted to sections of any involutive isotropic space $L$. This Lie bracket together with the projection $\pi_T :L \into TM$,  makes $L$ into a Lie algebroid and allows us to define a differential $d_L$ on $\gO^{\bullet}(L^*) = \Gamma(\wedge^{\bullet}L^*)$ making it into a DGA.
If $L$ is the $i$-eigenspace of a \gcs, then the natural pairing gives an isomorphism $L^* \cong \overline{L}$ and with this identification $(\gO^{\bullet}(\overline{L}),d_L)$ is a DGA.

If a \gcs\ has type zero over $M$, i.e., is of symplectic type, then  $\pi:L \stackrel{\cong}{\into} T_{\C}$ is an isomorphism and the Courant  bracket on $\Gamma(L)$ is mapped to the Lie bracket on $\Gamma(T_{\C})$. Therefore,  in this particular case, $(\gO^{\bullet}(\overline{L}),d_L)$ and  $(\gO_{\C}^{\bullet}(M),d)$ are isomorphic DGAs.

\subsection{Decomposition of forms}

A \gcs\ can also be
described using differential forms.  Recall that the exterior
algebra $\wedge^{\bullet}T^*$ carries a natural spin representation
for the metric bundle $T \oplus T^*$; the Clifford action of $X+\xi
\in T \oplus T^*$ on $ \rho \in \wedge^{\bullet}T^*$ is
$$ (X+\xi) \cdot \rho  = i_X \rho + \xi \wedge \rho.$$
The subspace $L\subset T_{\C} \oplus T_{\C}^*$ annihilating a spinor
$\rho \in \wedge^{\bullet}T_{\C}^*$ is always isotropic.  If $L$ is
maximal isotropic, then $\rho$ is called a {\it pure spinor} and it
must have the following algebraic form at every point:
\begin{equation}\label{eq:pure}
\rho =  e^{B+ i \omega}\wedge \gO,
\end{equation}
where $B$ and  $\omega$ are real 2-forms and $\gO$ is a decomposable
complex form.  Pure spinors annihilating the same space must be
equal up to rescaling, hence a maximal isotropic $L \subset T_{\C}
\oplus  T^*_{\C}$ may be uniquely described by a line subbundle
$U\subset \wedge^{\bullet}T_{\C}^*$.

For a complex manifold $U = \wedge^{n,0}T^*$ and for a symplectic manifold $U$ is generated by the globally defined closed form $e^{i \omega}$.  In general we have the following definition.

\begin{definition}
Given a \gcs\ \J, the line subbundle $U~\subset~\wedge^{\bullet}T^*_{\C}$ annihilating its $i$-eigenspace
is the {\it canonical bundle} of $\J$.
\end{definition}

Note that the condition $L \cap \overline{L} = \{0\}$ at the fiber of $E$ over $p \in M$ is equivalent
to the requirement that
\begin{equation}\label{eq:LcapLbar}
\gO \wedge \overline{\gO} \wedge \omega^{n-k} \neq 0
\end{equation}
for a generator $\rho = e^{B+ i \omega}\wedge \gO$ of $U$ at $p$, where $k = \deg(\gO)$ and $2n = \dim(M)$.

By letting $\wedge^{\bullet} \overline{L} \subset \Cliff(L \oplus \overline{L})$ act on the canonical line bundle we obtain a decomposition of the differential forms on $M^{2n}$: 
$$ \wedge^{\bullet}T^*_{\C}(M) = \oplus_{k=-n}^n U^k, \qquad \mbox{ where } U^{k} = \wedge^{n-k}\overline{L} \cdot \rho.$$
one can also describe the spaces $U^k$ as the $ik$-eigenspaces of $\J$ acting on forms.

Letting $\mc{U}^k = \Gamma(U^k)$ and $d_H = d + H\wedge$,  Courant integrability of the \gcs\ is equivalent to
$$d_H: \mc{U}^k \into  \mc{U}^{k+1} \oplus  \mc{U}^{k-1},$$ 
which allows us to define operators $\del: \mc{U}^k \into \mc{U}^{k+1}$ and $\delbar:\mc{U}^k \into \mc{U}^{k-1}$ by composing $d_H$ with the appropriate projections.

Given  a local section $\rho$ of the canonical bundle the operator $\delbar$ is related to $d_L$ by
$$ \delbar (\ga \cdot \rho)  = (d_L \ga) \cdot \rho + (-1)^{|\ga|}d_H\rho,$$
where $\ga \in \gO^{\bullet}(\overline{L})$ and $|\ga|$ is the degree of \ga. In the particular case when $(M,\J)$ has {\it holomorphically trivial canonical bundle}, i.e., there is a nonvanishing $d_H$-closed  global section $\rho$ of the canonical bundle, the above becomes
\begin{equation}\label{eq:dL and delbar}
\delbar (\ga \cdot \rho)  = (d_L \ga) \cdot \rho
\end{equation}
and hence $\rho$ furnishes an isomorphism of differential complexes. 

\section{Generalized K\"ahler manifolds}

In this section we introduce \gk\ manifolds. For these manifolds both $\Omega_{\C}^{\bullet}(M)$ and $(\Omega^{\bullet}(\overline{L}),d_L)$ admit a bigrading and, in certain conditions, some differential operators $\Omega_{\C}^{\bullet}(M)$ correspond to differential operators on $\Omega^{\bullet}(\overline{L})$. This correspondence was also used by Yi Li to study the moduli space of a \gks\ \cite{Li05} and is the key ingredient for our formality theorem. 

\begin{definition}
A {\it \gks} on a manifold $M^{2n}$ is a pair of commuting \gcss\ $\J_1$, $\J_2$  on $M$ such that
$$\IP{\J_1 \J_2 v, v} > 0 \qquad \mbox{ for } v \in T \oplus T^*\backslash\{0\}. $$
\end{definition}

Let $L_i$  be the $i$-eingenspace of $\J_i$. Since $\J_1$ and $\J_2$ communte, $\J_2$ furnishes a complex structure on $L_1$ with $i$-eigenspace $L_1\cap L_2 $. Using the fact that the natural pairing has signature $(n,n)$ and that $\IP{\J_1 \J_2 \cdot, \cdot }$ is positive definite one can show  $\dim(L_1) = 2\dim (L_1\cap L_2)$. Since $L_2$ is closed under the Courant bracket, we see that $L_1\cap L_2$ is closed under the bracket in the Lie algebroid $L_1$, and hence $\J_2|_{L_1}$ is an integrable complex structure on $L_1$. Using this complex structure we can decompose 
$$\wedge^{\bullet}\overline{L_1} = \oplus_{p,q}\wedge^{p,q}\overline{L_1}\qquad \mbox{ and }\qquad d_{L_1} = \del_{L_1} + \delbar_{L_1},$$
As in a complex manifold, the operators $\del_{L_1}$ and $\delbar_{L_1}$ are derivations, in the sense that they satisfy the Leibniz rule.

A \gk\ structure also gives a refinement of the deposition of forms into the spaces $U^k$. Since $\J_1$ and $\J_2$ commute one immediately obtains that the space of differential forms can be decomposed in terms of the eigenspaces of $\J_1$ and $\J_2$: $U^{p,q} = U_{\J_1}^p \cap U_{\J_2}^q$. This allows us to decompose $d_H$ further in 4 components
$$d_H: \mc{U}^{p,q} \into \mc{U}^{p+1,q+1} +  \mc{U}^{p+1,q-1}+\mc{U}^{p-1,q+1} +  \mc{U}^{p-1,q-1}.$$
In this case, the operator $\delbar$ for the \gcs\ $\J_1$  corresponds to the sum of the last two terms:
$$\delbar_1: \mc{U}^{p,q} \into \mc{U}^{p-1,q+1} +  \mc{U}^{p-1,q-1},$$
and we can define $\delta_+$ and $\delta_-$  as the projections of $\delbar_1$ into each of the components
$$\delta_+: \mc{U}^{p,q} \into \mc{U}^{p-1,q+1}\qquad \delta_-: \mc{U}^{p,q} \into \mc{U}^{p-1,q-1}.$$

By studying the Hodge theory of a \gk\ manifold, Gualtieri proved the following
\begin{theorem}\label{ddbar}
{\em (Gualtieri \cite{Gu04})} {\bf $\delta_+\delta_-$-lemma}. In a compact \gk\ manifold
$$\im \delta_+ \cap \Ker \delta_- = \im \delta_- \cap \Ker \delta_+ = \im(\delta_+\delta_-)  $$
\end{theorem}

If $\J_1$ has holomorphically trivial canonical bundle, then the correspondence between $\delbar_1$ and $d_{L_1}$ given in equation (\ref{eq:dL and delbar})  also furnishes  a correspondence between the operators  $\del_{L_1}$ and $\delbar_{L_1}$ on $\Omega^{\bullet}(\overline{L})$  and the operators
$\delta_+$ and $\delta_-$ on $\Omega_{\C}^{\bullet}(M)$. So, as a consequence of Theorem \ref{ddbar}, the operators $\del_{L_1}$ and $\delbar_{L_1}$ satisfy the $\del_{L_1}\delbar_{L_1}$-lemma and since they are derivations   the same argument from \cite{DGMS75} gives:
\begin{theorem}\label{gk formality}
If $(M,\J_1,\J_2)$ is a compact \gk\ manifold and $\J_1$ has holomophically trivial canonical bundle, then the DGA $(\Omega^{\bullet}(\overline{L_1}),d_{L_1})$ is formal.
\end{theorem}

In the case when $\J_1$ is a symplectic structure, then not only does it have a holomorphically trivial canonical bundle, but $(\Omega^{\bullet}(\overline{L_1}),d_{L_1})$ is isomorphic to $(\Omega_{\C}^{\bullet}(M),d)$. Therefore we have:

\begin{corollary}
If $(M,\J_1,\J_2)$ is a compact \gk\ manifold and $\J_1$ is a symplectic structure, then $M$ is formal.
\end{corollary}
This corollary generalizes the original theorem of formality of K\"ahler manifolds  \cite{DGMS75}.

\section{Nilpotent Lie algebras}

In this section we use Theorem \ref{gk formality} to prove that no nilpotent Lie algebra admits a \gk\ structure. Before we state the theorem we should stress that a \gcs\ on a Lie algebra $\frak{g}$ with closed 3-form $H \in \wedge^3 \frak{g}^*$ is just an integrable linear complex structure on $(\frak{g} \oplus \frak{g}^*, \Cour{~,~})$, orthogonal \wrt\ the natural pairing, where the Courant bracket is defined by
$$\Cour{X+\xi,Y+\eta}= [X,Y] + \mc{L}_X \eta -\mc{L}_Y\xi + i_Y i_X H,$$
and is a Lie bracket in this situation.

We also recall that a complex structure on a Lie algebra $\frak{g}$ is called {\it abelian} if its $i$-eigenspace, $\frak{g}^{1,0}$, is an abelian subalgebra of $\frak{g} \tensor \C$ \cite{BDM95,CFU02}.  By analogy, we say that a \gcs\ on $\frak{g}$ is {\it abelian} if the corresponding complex structure on $\frak{g}\oplus \frak{g}^*$ is abelian. Before we state our theorem on \gks s on nilpotent Lie algebras we need a little lemma:

\begin{lemma}\label{abelian gcss}
If a Lie algebra $\frak{g}$ admits an  abelian \gcs, then $\frak{g}$ is abelian.
\end{lemma}
\begin{proof}
Let $L$ be the $i$-eigenspce of an abelian \gcs\ on $\frak{g}$. Since $L$ is abelian, so is its projection over $\Gg\otimes \C$. Further, if $v \in \pi(L) \cap \pi(\overline{L})$ then $v$ is a central element in $\Gg_\C$. Indeed for such a $v$ there is a $\xi\in \Gg^*_\C$ such that $\J (v+\xi) \in \Gg^*_\C$ so, for $w \in \Gg_\C$
\begin{align*}
4[v,w] &= 4\pi(\Cour{v+\xi,w} )  \\ &= \pi(\Cour{v+\xi +i\J(v+\xi) + v+\xi - i\J(v+\xi), w + i\J w + w -i \J w}) \\
&= \pi(\Cour{v+\xi +i\J(v+\xi),w + i\J w} +\Cour{v+\xi +i\J(v+\xi),w-i\J w} +\\
&+\Cour{v+\xi -i\J(v+\xi),w + i\J w} +\Cour{v+\xi -i\J(v+\xi),w-i\J w})\\
& = \pi(\Cour{v+\xi +i\J(v+\xi),w-i\J w} +\Cour{v+\xi -i\J(v+\xi),w + i\J w})\\
&= \pi(\Cour{v+\xi -i\J(v+\xi),w-i\J w} +\Cour{v+\xi +i\J(v+\xi),w + i\J w})=0,
\end{align*}
where we have used that $L$ and $\overline{L}$ are abelian in the  fourth and in the last equalities and in the fifth equality we used that $\J(v+\xi) \in \Gg^*_\C$, hence the change of signs does not affect the projection of the bracket onto $\Gg_\C$.

If we let $e^{B+i \omega}\wedge \gO$, with $\gO = \theta^1 \wedge \cdots \wedge \theta^k$, be a generator for the canonical bundle of $\J$, then $\pi(L) \cap \pi(\overline{L})$ is the annihilator of $\gO \wedge \overline{\gO}$. Since $\theta^i \in L$ there are $\overline{\del_j} \in \overline{L}$ such that $\IP{\theta^i,\overline{\del_j}} = \delta^i_j$ and we can compute
$$\theta^i([\pi(\del_j),\pi(\overline{\del_k})]) = d\theta^i(\pi(\del_j),\pi(\overline{\del_k})) = \IP{\Cour{\theta^i,\del_j},\overline{\del_k}} = 0  $$
since $\theta^i, \del_j \in L$. Analogously we see that $[\pi(\del_j),\pi(\overline{\del_k})]$ also annihilates $\overline{\theta^i}$ and hence $[\pi(\del_j),\pi(\overline{\del_k})] \in \pi(L) \cap \pi(\overline{L})$, hence $\frak{g}$ is either abelian or 2-step nilpotent.

If $\frak{g}$ was 2-step nilpotent there would be an element $\xi \in \frak{g}^*$ with $d\xi \neq 0$. Since the only nonvanishing brackets are of the form $[\pi(\del_i),\pi(\delbar_j)]$ and $\xi$ is real, we see that there is a $\del_i$ for which $d\xi(\pi(\del_i),\pi(\delbar_i))= \xi([\pi(\del_i),\pi(\delbar_i)]) \neq 0$. Since all the $\theta^i$ are closed, we can further assume that $\xi = \J(v -B(v))$, for some $v \in \frak{g}$, therefore $v-B(v)-i\xi \in L$ and
\begin{align*}
0=\IP{\Cour{v-B(v)-i\xi, \del_i},\delbar_i} & = -id \xi(\pi(\del_i),\pi(\delbar_i)) + (H+dB)(v,\pi(\del_i),\pi(\delbar_i)).
\end{align*}
Observe that the first term is real and nonzero while the second is purely imaginary, hence the equation above can never hold and $\frak{g}$ is abelian.
\end{proof}

\begin{theorem}\label{nilpotent}
If a  nilpotent Lie algebra $\frak{g}$ admits a  \gk\ structure, then $\frak{g}$ is abelian.
\end{theorem}
\begin{proof}
According to  \cite{CG04a}, Theorem 3.1, every \gcs\ on a nilpotent Lie algebra $\Gg$  has holomorphically trivial canonical bundle. Further, for any closed $H \in \wedge^3 \Gg^*$, $\Gg\oplus \Gg^*$ with the Courant bracket  is again a nilpotent Lie algebra, hence the $i$-eigenspace, $L$, of any \gcss\ is a nilpotent Lie subalgebra of $(\Gg\oplus \Gg^*)\tensor \C$. According to Lemma \ref{abelian gcss}, if $\frak{g}$ has nontrivial bracket, then $L$ has a nontrivial bracket. Then, Example \ref{ex:nilpotent lie algebras} shows that $(\wedge^{\bullet}\overline{L},d_L)$ is not formal and hence, by Theorem \ref{gk formality}, can not be part of a \gk\ pair.
\end{proof}


\end{document}